\documentclass[12pt]{article}
\usepackage{cite}
\usepackage{graphicx}
\usepackage[cmex10]{amsmath}
\usepackage{amssymb}
\usepackage{arydshln}
\usepackage{multirow}
\usepackage{caption}
\usepackage[hidelinks]{hyperref}
\newcommand{\R}{\mathbb{R}}
\newcommand{\Z}{\mathbb{Z}}
\newtheorem{thm}{Theorem}
\newtheorem{cor}{Corollary}
\newtheorem{rem}{Remark}
\newtheorem{lem}{Lemma}
\newtheorem{prop}{Proposition}

\begin{document}

\title{Predictor-based networked control \\under uncertain transmission delays}

\author{Anton~Selivanov \and
        Emilia~Fridman
\thanks{A. Selivanov ({\tt\small antonselivanov@gmail.com}) and E. Fridman ({\tt\small emilia@eng.tau.ac.il}) are with School of Electrical Engineering, Tel Aviv University, Israel.}% <-this % stops a space
\thanks{This work was supported by Israel Science Foundation (grant No. 1128/14) and has been published in \cite{Selivanov2016e}.}%
}
\renewcommand\footnotemark{}
\maketitle

\begin{abstract}
  We consider state-feedback predictor-based control of networked control systems with large time-varying communication delays. We show that even a small controller-to-actuators delay uncertainty may lead to a non-small residual error in a networked control system and reveal how to analyze such systems. Then we design an event-triggered predictor-based controller with sampled measurements and demonstrate that, depending on the delay uncertainty, one should choose various predictor models to reduce the error due to triggering. For the systems with a network only from a controller to actuators, we take advantage of the continuously available measurements by using a continuous-time predictor and employing a recently proposed switching approach to event-triggered control. By an example of an inverted pendulum on a cart we demonstrate that the proposed approach is extremely efficient when the uncertain time-varying network-induced delays are too large for the system to be stabilizable without a predictor. 
\end{abstract}

\section{Introduction}
In networked control systems (NCSs), which are comprised of sensors, controllers, and actuators connected through a communication medium, transmitted signals are sampled in time and are subject to time-delays. Most existing papers on NCSs study robust stability with respect to small communication delays (see, e.g., \cite{Antsaklis2004,Fridman2004,Gao2008,Liu2012b}). To compensate large transport delays, predictor-based approach can be employed. So far this was done only for sampled-data control with \textit{known constant delays} \cite{Karafyllis2012a,Mazenc2013}. In this paper we develop predictor-based sampled-data control for \textit{unknown time-varying delays}.

There are several works that study robustness (w.r.t.~delay uncertainty) of a predictor-based \textit{continuous-time} controller \cite{Yue2005,Bekiaris-liberis2013,Karafyllis2013,Li2014}. In these works the residual error that appears due to delay uncertainty can be made arbitrary small by reducing the upper bound of the uncertainty. However, this is not true for \textit{sampled-data systems}, where an arbitrary small delay uncertainty may lead to a non-vanishing error (because the terms that appear in the residual error may belong to different sampling intervals).

In this work we study an NCS with two networks: from sensors to a controller and from the controller to actuators. Both networks introduce large time-varying delays. We assume that the messages sent from the sensors are time stamped \cite{Zhang2001}. This allows to calculate the sensors-to-controller delay. The controller-to-actuators delay is assumed to be unknown but belongs to a known delay interval. We use a state-feedback predictor, which is calculated on the controller side, to partially compensate both delays. By extending the time-delay modelling of NCSs \cite{Fridman2004,Gao2008,Fridman2014}, we present the system in a form suitable for analysis. Using a proper Lyapunov-Krasovskii functional, we derive LMI-based conditions for the stability analysis and design that guarantee the desired decay rate of convergence.

As the next step we introduce an event-triggering mechanism \cite{Tabuada2007,Heemels2012} into predictor-based networked control. The event-triggering condition is checked on a controller side and allows to reduce the amount of control signals sent through a controller-to-actuators network. We demonstrate that it is reasonable to choose different predictor models for a zero and non-zero controller-to-actuators delay uncertainty. Finally, we consider predictor-based event-triggered control with continuous-time measurements and sampled control signals sent through a controller-to-actuators network. Such systems naturally appear when a visually observed vehicle is controlled through a wireless network. To take advantage of the continuously available measurements, we use a continuous-time predictor \cite{Mazenc2013,Kwon1980,Artstein1982} and a recently proposed switching approach to event-triggered control~\cite{Selivanov2016c}. 

By an example of an inverted pendulum on a cart we demonstrate that the proposed approach is extremely efficient when the uncertain time-varying network-induced delays are too large for the system to be stabilizable without a predictor. Moreover, the considered event-triggering mechanism allows to significantly reduce the network workload. 

\section{Networked control employing predictor}\label{sec:pred}
\begin{figure}
	\centering
	\includegraphics[width=.6\linewidth]{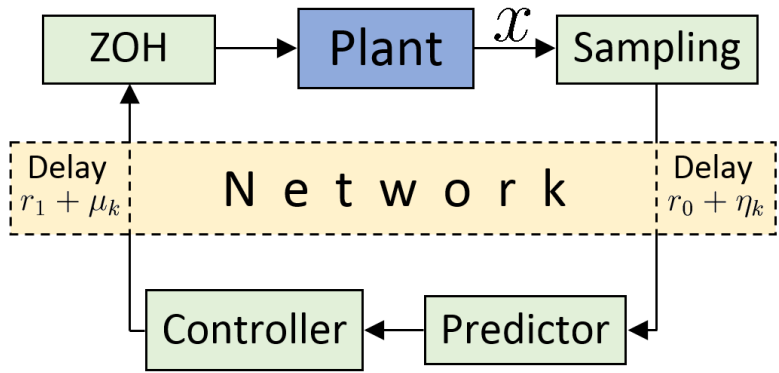}
	\caption{NCS with a predictor}\label{fig:scheme1}
\end{figure}
\begin{figure}
	\centering
	\includegraphics[width=.6\linewidth]{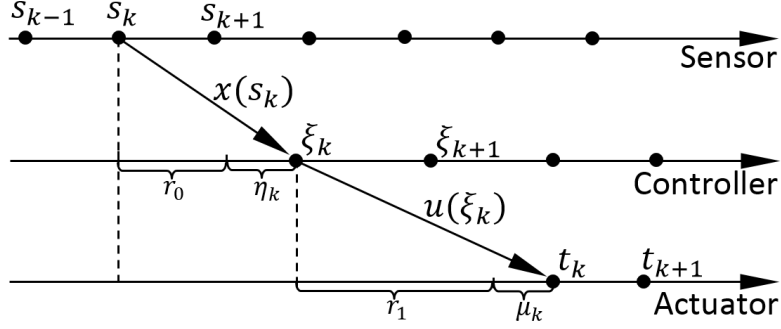}
	\caption{Time-delays and updating times}\label{updatingTimes}
\end{figure}

Consider the linear system
\begin{equation}\label{eq:inSys}
	\dot x(t)=Ax(t)+Bu(t),\quad t\ge0
\end{equation}
with the state $x\in\R^n$, control input $u\in\R^m$, and constant matrices $A$, $B$ of appropriate dimensions for which there exists $K\in\R^{m\times n}$ such that $A+BK$ is a Hurwitz matrix. Let $\{s_k\}$ be sampling instants such that
\begin{equation*}
	\begin{array}{c}
		0=s_0<s_1<s_2<\ldots,\qquad\lim_{k\to\infty}s_k=\infty,\qquad s_{k+1}-s_k\leq h.
	\end{array}
\end{equation*}
At each sampling time $s_k$ the state $x(s_k)$ is transmitted to a controller, where a control signal is calculated and transmitted to actuators (see Fig.~\ref{fig:scheme1}). We assume that the controller and the actuators are event-driven (update their outputs as soon as they receive new data). Both state and control signals are subject to network-induced delays $r_0+\eta_k$ and $r_1+\mu_k$, respectively. Thus, the controller updating times are $\xi_k=s_k+r_0+\eta_k$ and the actuators updating times are $t_k=\xi_k+r_1+\mu_k$, where $k\in\Z_+$, $\Z_+=\{0,1,2,\ldots\}$ (see Fig.~\ref{updatingTimes}). Here $r_0$ and $r_1$ are known constant transport delays, $\eta_k$ and $\mu_k$ are time-varying delays such that
\begin{equation}\label{etamu}
	0\leq\eta_k\le\eta_M,\qquad0\leq\mu_k\le\mu_M,\qquad\xi_k\leq\xi_{k+1},\qquad t_k\leq t_{k+1}.
\end{equation}
We assume that the sensors and controller clocks are synchronized and together with $x(s_k)$ the time stamp $s_k$ is transmitted so that the value of $\eta_k=\xi_k-s_k-r_0$ can be calculated on the controller side at time $\xi_k$. Delay uncertainty $\mu_k$ is assumed to be unknown. Note that we do not require $\eta_k+\mu_k$ to be less than the sampling interval but the sequences $\{\xi_k\}$ and $\{t_k\}$ of updating times should be increasing.

Define $u(\xi)=0$ for $\xi<\xi_0$. Then \eqref{eq:inSys} transforms to
\begin{equation}\label{eq:sampSys}
	\begin{aligned}
		\dot x(t)&=Ax(t),&&t\in[0,t_0), \\
		\dot x(t)&=Ax(t)+Bu(\xi_k),&&t\in[t_k,t_{k+1}),\quad k\in\Z_+.
	\end{aligned}
\end{equation}
To construct a predictor-based controller for \eqref{eq:sampSys}, define
\begin{equation}\label{vdef}
	v(\xi)\triangleq\left\{\begin{aligned}
		&0,&&\xi<\xi_0, \\
		&u(\xi_k),&&\xi\in[\xi_k,\xi_{k+1}),\quad k\in\Z_+
	\end{aligned}\right.
\end{equation}
and consider the change of variable \cite{Kwon1980,Artstein1982}
\begin{equation}\label{pred}
		z(t)\triangleq e^{A(r_0+r_1)}x(t)+\int_{t-r_1}^{t+r_0}e^{A(t+r_0-\theta)}Bv(\theta)\,d\theta,
\end{equation}
where $t\ge0$. We set $z(t)=0$ for $t<0$. If $\mu_M=0$, i.e. controller-to-actuators delay is constant, \eqref{vdef}, \eqref{pred} is the state prediction, namely, $z(t)=x(t+r_0+r_1)$. If $\mu_k\not\equiv0$ to obtain the precise state prediction one needs to integrate \eqref{eq:sampSys}, where $t_k=\xi_k+r_1+\mu_k$ depends on $\mu_k$. Since $\mu_k$ is unknown, we use the prediction \eqref{vdef}, \eqref{pred} that is imprecise for $\mu_k\not\equiv0$. By substituting \eqref{eq:sampSys} for $\dot x(t)$ we obtain
\begin{equation}\label{eq:vzSys}
	\begin{aligned}
		\dot z(t)&=Az(t)+Bv(t+r_0)-e^{A(r_0+r_1)}Bv(t-r_1),&&t\in[0,t_0), \\
		\dot z(t)&=Az(t)+Bv(t+r_0)+e^{A(r_0+r_1)}B\left[u(\xi_k)-v(t-r_1)\right],&&t\in[t_k,t_{k+1}),\quad k\in\Z_+.
	\end{aligned}
\end{equation}
Consider the following control law
\begin{equation}\label{cont}
		u(\xi_k)\triangleq Kz(s_k)=K\left[e^{A(r_0+r_1)}x(s_k)
		+\int_{\xi_k-\eta_k-r_0-r_1}^{\xi_k-\eta_k}e^{A(\xi_k-\eta_k-\theta)}B v(\theta)\,d\theta\right],\quad k\in\Z_+.
\end{equation}
Since $\eta_k$ is available to the controller at time $\xi_k$, the control signal \eqref{cont} can be calculated. Moreover, no numerical difficulties arise while calculating the integral term in \eqref{cont} with a piecewise constant $v(\theta)$ given by \eqref{vdef}.

We analyze \eqref{vdef}--\eqref{cont} using the time-delay approach to NCSs \cite{Fridman2004,Gao2008,Fridman2014}. According to \eqref{vdef}, \eqref{cont}, $v(t+r_0)=Kz(s_k)$ whenever $t+r_0\in[\xi_k,\xi_{k+1})$, that is, when $t\in[\xi_k-r_0,\xi_{k+1}-r_0)$. If $t<\xi_0-r_0$ then $v(t+r_0)=0=Kz(t-\eta_0)$. Therefore,
\begin{equation}\label{eq:9}
	v(t+r_0)=Kz(t-\tau(t)),\quad t\in\R,
\end{equation}
where
\begin{equation*}
	\tau(t)=\left\{\begin{aligned}
		&\eta_0,&&t<\xi_0-r_0, \\
		&t-s_k,&&t\in[\xi_k-r_0,\xi_{k+1}-r_0),\quad k\in\Z_+.
	\end{aligned}\right.
\end{equation*}
Note that for $t\ge\xi_0-r_0$
\begin{equation*}
	0\leq\tau(t)\leq\max_k\{(s_{k+1}+r_0+\eta_{k+1})-r_0-s_k\}\leq h+\eta_M.
\end{equation*}
By similar reasoning we obtain
\begin{equation}\label{eq:cll}
\dot z(t)=Az(t)+BKz(t-\tau(t))+e^{A(r_0+r_1)}BK[z(t-\tau_2(t))-z(t-\tau_1(t))],\quad t\ge0,
\end{equation}
with
\begin{equation}\label{inCond}
	z(0)=e^{A(r_0+r_1)}x(0), \quad z(t)=0\text{ for }t<0,
\end{equation}
\begin{equation}\label{delays}
	\begin{aligned}
		\tau(t)&\triangleq\left\{\begin{aligned}
			&\eta_0,&&t<\xi_0-r_0, \\
			&t-s_k,&&t\in[\xi_k-r_0,\xi_{k+1}-r_0), k\in\Z_+,
		\end{aligned}\right.\\
		\tau_1(t)&\triangleq \left\{\begin{aligned}
			&r_1+r_0+\eta_0,&&t\in[0,t_0-\mu_0),\\
			&t-s_k,&&t\in[t_k\!-\!\mu_k,t_{k+1}\!-\!\mu_{k+1}), k\!\in\!\Z_+,
		\end{aligned}\right.\\
		\tau_2(t)&\triangleq\left\{\begin{aligned}
			&r_0+r_1+\eta_0+\mu_0,&&t\in[0,t_0), \\
			&t-s_k,&&t\in[t_k,t_{k+1}), k\in\Z_+,
		\end{aligned}\right.
	\end{aligned}
\end{equation}
\begin{equation*}
	\begin{array}{c}
		0\leq\tau(t)\leq\bar\tau\triangleq h+\eta_M, \\
		r_0+r_1\leq\tau_1(t)\leq\tau_2(t)\leq\tau_M\triangleq r_0+r_1+h+\eta_M+\mu_M.
	\end{array}
\end{equation*}

\begin{rem}
	If $\xi_k=\xi_{k+1}$ then $\tau(t)=t-s_{k-1}$ for $t\in[\xi_{k-1}-r_0,\xi_{k+1}-r_0)$ and it may seem that the bound $\tau(t)\le h+\eta_M$ can be violated. This is not the case, since $\xi_k=\xi_{k+1}$ implies $s_k+r_0+\eta_k=s_{k+1}+r_0+\eta_{k+1}$, that is, $\eta_{k+1}\le \eta_k-h$. Therefore, for $t\in[\xi_{k-1}-r_0,\xi_{k+1}-r_0)$
	\begin{multline*}
		\tau(t)\le\xi_{k+1}-r_0-s_{k-1}=s_{k+1}+r_0+\eta_{k+1}-r_0-s_{k-1}\\
		\le (s_{k+1}-s_{k-1})+(\eta_k-h)\le 2h+\eta_k-h=\eta_k+h.
	\end{multline*}
	Similar explanation is valid for $\xi_k=\xi_{k+1}=\cdots=\xi_{k+d}$ and $t_k=t_{k+1}=\cdots=t_{k+d}$.
\end{rem}

\begin{rem}
	If $\mu_k\!\equiv\!0$ then $\tau_1(t)\!=\tau_2(t)$ and \eqref{eq:cll} simplifies to
	\begin{equation}\label{eq:noMu}
		\dot z(t)=Az(t)+BKz(t-\tau(t)),\quad t\geq 0.
	\end{equation}
	The system \eqref{eq:noMu} is independent of $r_0$ and $r_1$. Therefore, the stability conditions for \eqref{eq:noMu} are independent of $r_0$ and $r_1$: these delays are compensated by the predictor \eqref{vdef}, \eqref{pred}. For $\mu_k\not\equiv0$ the system \eqref{eq:cll} contains the residual error that appears due to impreciseness of the predictor \eqref{vdef},~\eqref{pred}.
\end{rem}

\begin{rem}
	While studying robustness of a predictor for the time-delay $r+\mu(t)$ with the uncertainty $\mu(t)\leq\mu_M$, usually the residual $e^{Ar}BK[z(t-r-\mu(t))-z(t-r)]$ appears in the closed-loop system~\cite{Karafyllis2013,Fridman2014}. Since $\dot z$ is generally proved to be bounded, even for unstable $A$ and large $r$ by reducing $\mu_M$ one can retain this error small enough to preserve the stability. In a word, $r$ can be made arbitrary large by decreasing $\mu_M$. This does not hold for sampled-data systems: for arbitrary small $\mu_k>0$ when $t\in[t_k-\mu_k,t_k)$ the arguments of $z(t-\tau_1(t))$ and $z(t-\tau_2(t))$ belong to different sampling intervals, namely, $(t-\tau_1(t))-(t-\tau_2(t))=s_k-s_{k-1}$ (if $t_k-\mu_k>t_{k-1}$, $k\ge1$). Therefore, smallness of the residual in \eqref{eq:cll} for large $r=r_0+r_1$ can be guaranteed only by reducing $\mu_M$ together with the maximum sampling interval $h$.
\end{rem}

Stability conditions for the systems \eqref{eq:cll} and \eqref{eq:noMu} follow from Theorem~\ref{th:1} and Proposition~\ref{prop:1_2} of the next section. 

\section{Event-triggering with sampled measurements}\label{sec:evtr}
\begin{figure}
	\centering
	\includegraphics[width=.6\linewidth]{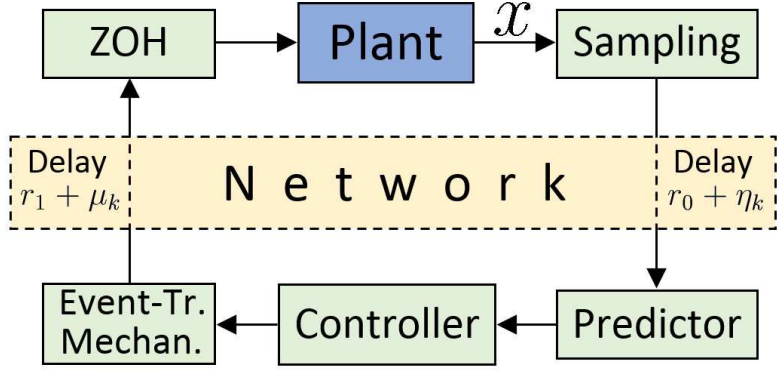}
	\caption{NCS with a predictor and event-triggering mechanism}\label{fig:scheme2}
\end{figure}
To reduce the workload of a controller-to-actuators network, we incorporate an event-triggering mechanism (see \cite{Tabuada2007}). The idea is to send only those control signals $u(\xi_k)$ which relative change is greater than a constant $\sigma\in[0,1)$ (see Fig.~\ref{fig:scheme2}), namely, that violate the following event-triggering rule
\begin{equation}\label{EvTrRe}
\left(\hat u(\xi_{k-1})-u(\xi_k)\right)^T\Omega\left(\hat u(\xi_{k-1})-u(\xi_k)\right)\leq\sigma u^T(\xi_k)\Omega u(\xi_k),
\end{equation}
where a matrix $\Omega\ge0$ and a scalar $\sigma\ge0$ are event-triggering parameters and $\hat u(\xi_{k-1})$ is the last sent control value before the time instant $\xi_k$:
\begin{equation}\label{EvTrCont}
\hat u(\xi_k)=\left\{
\begin{aligned}
&\hat u(\xi_{k-1}),&& \text{if \eqref{EvTrRe} is true}, \\
&u(\xi_k),&& \text{otherwise,}
\end{aligned}
\right.
\end{equation}
with $\hat u(\xi_{-1})=0$. Note that the sensor sends measurements at time instants $s_k$ (such that $s_{k+1}-s_k\le h$) independent of the event-triggering events. Then \eqref{eq:sampSys} takes the form
\begin{equation}\label{eq:startSys}
\begin{aligned}
\dot x(t)&=Ax(t),&&t\in[0,t_0),\\
\dot x(t)&=Ax(t)+B\hat u(\xi_k),&&t\in[t_k,t_{k+1}),\quad k\in\Z_+.
\end{aligned}
\end{equation}
Consider the change of variable \eqref{pred} with $v(\theta)$ to be defined. By substituting \eqref{eq:startSys} for $\dot x(t)$, we obtain
\begin{equation}\label{eq:vzzSys}
\hspace{0em}\begin{aligned}
\dot z(t)&=Az(t)+Bv(t+r_0)-e^{A(r_0+r_1)}Bv(t-r_1),\,t\in[0,t_0), \\
\dot z(t)&=Az(t)+Bv(t+r_0)+e^{A(r_0+r_1)}B\left[\hat u(\xi_k)-v(t-r_1)\right],\quad t\in[t_k,t_{k+1}),\quad k\in\Z_+.
\end{aligned}
\end{equation}
We now show that for $\mu_M=0$ and $\mu_M>0$ one should pick different functions $v(\theta)$ in the predictor \eqref{pred}.

1. Let $\mu_M=0$. To cancel the last term in \eqref{eq:vzzSys} we take $v(t-r_1)=\hat u(\xi_k)$ for $t\in[t_k,t_{k+1})$ or, equivalently,
\begin{equation}\label{vmu0}
v(\xi)\triangleq\left\{\begin{aligned}
&0,&&\xi<\xi_0,\\
&\hat u(\xi_k),&&\xi\in[\xi_k,\xi_{k+1}),\quad k\in\Z_+.
\end{aligned}\right.
\end{equation}
Then \eqref{pred}, \eqref{vmu0} is the state prediction for the system~\eqref{eq:startSys}, i.e. $z(t)=x(t+r_0+r_1)$. The system \eqref{eq:vzzSys} takes the form
\begin{equation*}
\dot z(t)=Az(t)+Bv(t+r_0),\quad t\ge0.
\end{equation*}
Let us define
\begin{equation*}
e_0(t)\triangleq\left\{\begin{aligned}
&0,&&t<\xi_0,\\
&\hat u(\xi_k)-u(\xi_k),&&t\in[\xi_k,\xi_{k+1}),\quad k\in\Z_+.
\end{aligned}\right.
\end{equation*}
Then for $t\in[\xi_k-r_0,\xi_{k+1}-r_0)$ we have
\begin{equation*}
v(t+r_0)=\hat u(\xi_k)=u(\xi_k)+e_0(t+r_0)=Kz(s_k)+e_0(t+r_0)=Kz(t-\tau(t))+e_0(t+r_0)
\end{equation*}
with $\tau(t)$ defined in \eqref{delays}. For $t<\xi_0-r_0$, $v(t+r_0)=0=Kz(t-\eta_0)+e_0(t+r_0)$. Therefore,
\begin{equation}\label{eq:1}
\dot z(t)=Az(t)+BKz(t-\tau(t))+Be_0(t+r_0),\quad t\ge0
\end{equation}
with \eqref{inCond}, and \eqref{EvTrRe}, \eqref{EvTrCont} yield
\begin{equation*}
0\le\sigma z^T(t-\tau(t))K^T\Omega Kz(t-\tau(t))-e_0^T(t+r_0)\Omega e_0(t+r_0)
\end{equation*}
for $t\ge0$. It may seem that \eqref{eq:1} depends on the future, since $e_0(t+r_0)$ enters the system. This is not the case, since $e_0(\xi)$ for $\xi\in[\xi_k,\xi_{k+1})$ is fully defined by $z(s)$ with $s\le s_k=\xi_k-r_0-\eta_k$.

2. Let $\mu_k\not\equiv0$. Then the last term in \eqref{eq:vzzSys} cannot be canceled, since this would require to take $v(\xi)=\hat u(\xi_k)$ for $\xi\in[\xi_k+\mu_k,\xi_{k+1}+\mu_{k+1})$ with unknown $\mu_k$. If one defines $v(\xi)$ as in \eqref{vmu0} and uses $v(\xi)=\hat u(\xi_k)=u(\xi_k)+e_0(\xi)$, the functions $v(t+r_0)$, $v(t-r_1)$, $\hat u(\xi_k)$ present in \eqref{eq:vzzSys} will introduce three errors due to triggering $e_0$ with different arguments. To avoid additional triggering errors, we don't include them into the definition of $v(\xi)$, namely, we use \eqref{vdef} where $v(\xi)=u(\xi_k)$ or zero. Let us define
\begin{equation*}
e_1(t)\triangleq\left\{\begin{aligned}
&0,&&t<t_0,\\
&\hat u(\xi_k)-u(\xi_k),&&t\in[t_k,t_{k+1}),\quad k\in\Z_+.
\end{aligned}\right.
\end{equation*}
Then we have
\begin{equation*}
\begin{aligned}
0&=Kz(t-\tau_2(t))+e_1(t),\quad t\in[0,t_0), \\
\hat u(\xi_k)&=u(\xi_k)+e_1(t)=Kz(s_k)+e_1(t)=Kz(t-\tau_2(t))+e_1(t),\quad t\in[t_k,t_{k+1}),\quad k\in\Z_+
\end{aligned}
\end{equation*}
with $\tau_2(t)$ defined in \eqref{delays}. By arguments similar to those from Section~\ref{sec:pred} we obtain
\begin{equation}\label{eq:zSys}
\dot z(t)=Az(t)+BKz(t-\tau(t))+e^{A(r_0+r_1)}Be_1(t)+e^{A(r_0+r_1)}BK[z(t-\tau_2(t))-z(t-\tau_1(t))],\,t\ge0,
\end{equation}
with \eqref{inCond}, where $\tau$, $\tau_1$, $\tau_2$ are defined in \eqref{delays} and due to \eqref{EvTrRe}, \eqref{EvTrCont} for $t\ge0$
\begin{equation}\label{EvTrCompensation}
0\leq\sigma z^T(t-\tau_2(t))K^T\Omega Kz(t-\tau_2(t))-e_1^T(t)\Omega e_1(t).
\end{equation}
\begin{rem}
	Note that for $\mu_M=0$ \eqref{eq:zSys} transforms to
	\begin{equation*}
	\dot z(t)=Az(t)+BKz(t-\tau(t))+e^{A(r_0+r_1)}Be_1(t).
	\end{equation*}
	Since the triggering error $e_1(t)$ is multiplied by $e^{A(r_0+r_1)}$, to guarantee the stability of the system for unstable $A$ and large $r_0+r_1$ one needs to retain $e_1(t)$ small enough. This problem doesn't appear in the system \eqref{eq:1} for which the stability conditions are independent of $r_0$ and $r_1$ (see Proposition~\ref{prop:1_2}).
\end{rem}

\begin{figure}
	\centering
	\includegraphics[width=.8\linewidth]{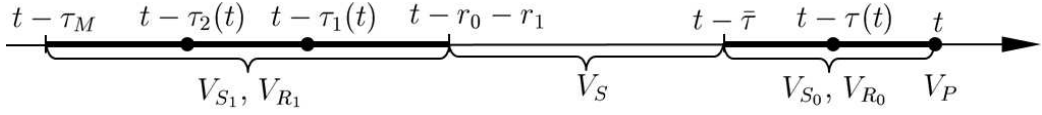}
	\caption{Lyapunov-Krasovskii functional}\label{Fig:LKfun}
\end{figure}
To avoid some technical complications, we assume that $\bar\tau=h+\eta_M\leq r_0+r_1$. The stability conditions are derived using Lyapunov-Krasovskii functional (see Fig.~\ref{Fig:LKfun})
\begin{equation*}
V=V_P+V_{S_0}+V_{R_0}+V_S+V_{S_1}+V_{R_1},
\end{equation*}
where
\begin{equation*}
\begin{array}{rl}
V_P&=z^T(t)Pz(t),\quad P>0,\\
V_{S_0}&=\int_{t-\bar\tau}^te^{2\alpha(s-t)}z^T(s)S_0z(s)\,ds,\quad S_0\ge0,\\
V_{R_0}&=\bar\tau\int_{-\bar\tau}^0\int_{t+\theta}^te^{2\alpha(s-t)}\dot z^T(s)R_0\dot z(s)\,ds\,d\theta,\quad R_0\ge0,\\
V_S&=\int_{t-r_0-r_1}^{t-\bar\tau}e^{2\alpha(s-t)}z^T(s)Sz(s)\,ds,\quad S\ge0,\\
V_{S_1}&=\int_{t-\tau_M}^{t-r_0-r_1}e^{2\alpha(s-t)}z^T(s)S_1z(s)\,ds,\quad S_1\ge0,\\
V_{R_1}&=(\tau_M-r_0-r_1)\int_{-\tau_M}^{-r_0-r_1}\int_{t+\theta}^te^{2\alpha(s-t)}\dot z^T(s)R_1\dot z(s)\,ds\,d\theta,\quad R_1\ge0.
\end{array}
\end{equation*}
Note that the delayed arguments of $z$ in \eqref{eq:zSys} belong to two bold regions in Fig.~\ref{Fig:LKfun}. To analyze these regions, we use standard delay-dependent terms in $V$ (see, e.g., \cite{Fridman2014}). To allow for large transport delays $r_0$ and $r_1$, we use only delay-independent term $V_S$ for the interval $[t-r_0-r_1,t-\bar\tau]$.
\begin{lem}\label{lem:1}
	For given $\mu_M\ge0$, $\eta_M\ge0$, and $\alpha>0$ let there exist an $n\times n$ matrix $P>0$, $n\times n$ non-negative matrices $S$, $S_0$, $S_1$, $R_0$, $R_1$, an $m\times m$ matrix $\Omega\ge0$, and $n\times n$ matrices $P_2$, $P_3$, $G_i$ ($i=0,\ldots,3$) such that
	\begin{equation*}
	\Phi\leq0,\quad\left[\begin{smallmatrix}
	R_0 & G_0 \\
	* & R_0
	\end{smallmatrix}\right]\ge0,\quad
	\left[\begin{smallmatrix}
	R_1 & G_i \\
	* & R_1
	\end{smallmatrix}\right]\ge0,\quad i=1,2,3,
	\end{equation*}
	where $\Phi=\{\Phi_{ij}\}$ is the symmetric matrix composed from
	\begin{equation*}
	\begin{aligned}
	\Phi_{11}&=2\alpha P+S_0-\bar\rho R_0+P_2^TA+A^TP_2,&& \Phi_{44}=\bar\rho(S-S_0-R_0),\\
	\Phi_{12}&=P-P_2^T+A^TP_3,&&\Phi_{55}=e^{-2\alpha(r_0+r_1)}(S_1-S)-\rho_MR_1,\\
	\Phi_{13}&=\bar\rho(R_0-G_0)+P_2^TBK,&&\Phi_{56}=\rho_M(R_1-G_1),\\
	\Phi_{14}&=\bar\rho G_0,&&\Phi_{57}=\rho_M(G_1-G_2),\\
	\Phi_{19}&=P_2^Te^{A(r_0+r_1)}B,&&\Phi_{58}=\rho_MG_2,\\
	\Phi_{17}&=-\Phi_{16}=\Phi_{19}K,&&\Phi_{66}=-\Phi_{56}-\Phi_{56}^T, \\
	\Phi_{22}&=\bar\tau^2R_0+(\tau_M-r_0-r_1)^2R_1-P_3-P_3^T,&& \Phi_{67}=\rho_M(R_1-G_1+G_2-G_3),\\
	\Phi_{23}&=P_3^TBK,&&\Phi_{68}=\rho_M(G_3-G_2), \\
	\Phi_{29}&=P_3^Te^{A(r_0+r_1)}B,&&\Phi_{78}=\rho_M(R_1-G_3),\\
	\Phi_{27}&=-\Phi_{26}=\Phi_{29}K,&&\Phi_{77}=-\Phi_{78}-\Phi_{78}^T+\sigma K^T\Omega K,\\
	\Phi_{34}&=\bar\rho(R_0-G_0),&&\Phi_{88}=-\rho_M(S_1+R_1),\\
	\Phi_{33}&=-\Phi_{34}-\Phi_{34}^T,&&\Phi_{99}=-\Omega,\\
	\end{aligned}
	\end{equation*}
	$\bar\rho=e^{-2\alpha\bar\tau}$, $\rho_M=e^{-2\alpha\tau_M}$, other blocks are zero matrices. Then the system \eqref{inCond}, \eqref{eq:zSys} is exponentially stable with a decay rate~$\alpha$, i.e. for some $M>0$ solutions of the system satisfy
	\begin{equation}\label{zdec}
	|z(t)|\leq Me^{-\alpha t}|z(0)|, \quad t\ge 0.
	\end{equation}
\end{lem}
{\em Proof} is given in Appendix~\ref{app:lem:1}.

\begin{thm}[Sampled event-triggering]\label{th:1}
	Under the conditions of Lemma~\ref{lem:1} the system \eqref{cont}, \eqref{EvTrRe}--\eqref{eq:startSys} with $v(\theta)$ given by \eqref{vdef} is exponentially stable with a decay rate $\alpha$, i.e. for some $M>0$ solutions of the system satisfy
	\begin{equation}\label{xdec}
	|x(t)|\le Me^{-\alpha t}|x(0)|.
	\end{equation}
	
\end{thm}
{\em Proof} is given in Appendix~\ref{app:th:1}.

\begin{rem}
	If $A+BK$ is Hurwitz and $\alpha=\tau_M=0$ the LMIs of Lemma~\ref{lem:1} are always feasible by the standard arguments for delay-dependent conditions \cite{Fridman2014}. That is, LMIs of Lemma~\ref{lem:1} establish relation between the decay rate, sampling period, and time-delays that preserve exponential stability of the system \eqref{vdef}, \eqref{cont}, \eqref{EvTrRe}--\eqref{eq:startSys}.
\end{rem}

\begin{cor}\label{cor:2}
	If conditions of Lemma~\ref{lem:1} are satisfied with $\sigma=0$, $\Omega>0$, the system \eqref{eq:sampSys} under the control law \eqref{cont} with $v(\theta)$ given by \eqref{vdef} is exponentially stable with a decay rate $\alpha$.
\end{cor}
{\em Proof.} For $\sigma=0$, $\Omega>0$ event-triggering mechanism \eqref{EvTrRe}, \eqref{EvTrCont} implies $\hat u(\xi_k)=u(\xi_k)$ and $e_1(t)\equiv0$, therefore, \eqref{eq:zSys} coincides with \eqref{eq:cll}. Then under conditions of Lemma~\ref{lem:1} \eqref{eq:cll} is exponentially stable. This implies exponential stability of \eqref{eq:sampSys}, \eqref{vdef}, \eqref{cont} due to the change of variable \eqref{vdef}, \eqref{pred}.

\hfill$\square$\\
For the case of $\mu_M=0$ the next proposition gives stability conditions independent of $r_0$ and~$r_1$.
\begin{prop}\label{prop:1_2}
	For $\mu_M=0$ and given $\eta_M\ge0$, $\alpha>0$, if there exist an $n\times n$ matrix $P>0$, $n\times n$ non-negative matrices $S$, $R$, an $m\times m$ matrix $\Omega\ge0$, and $n\times n$ matrices $P_2$, $P_3$, $G$ such that
	\begin{equation*}
	\Psi\leq0,\quad\left[\begin{smallmatrix}
	R & G \\
	* & R
	\end{smallmatrix}\right]\ge0,
	\end{equation*}
	where $\Psi=\{\Psi_{ij}\}$ is the symmetric matrix composed from
	\begin{align*}
	\Psi_{11}&=2\alpha P+S-\bar\rho R+P_2^TA+A^TP_2,&\Psi_{25}&=P_3^TB,\\
	\Psi_{12}&=P-P_2^T+A^TP_3,&\Psi_{23}&=\Psi_{25}K,\\
	\Psi_{13}&=\bar\rho(R-G)+P_2^TBK,&\Psi_{34}&=\bar\rho(R-G),\\
	\Psi_{14}&=\bar\rho G,&\Psi_{33}&=-\Psi_{34}-\Psi_{34}^T+\sigma K^T\Omega K,\\
	\Psi_{15}&=P_2^TB,&\Psi_{44}&=-\bar\rho(S+R),\\
	\Psi_{22}&=\bar\tau^2R-P_3-P_3^T,&\Psi_{55}&=-\Omega,\\
	\end{align*}
	$\bar\rho=e^{-2\alpha\bar\tau}$, other blocks are zero matrices, then \eqref{cont}, \eqref{EvTrRe}--\eqref{eq:startSys} with $v(\theta)$ given by \eqref{vmu0} is exponentially stable with a decay rate~$\alpha$.
\end{prop}
{\em Proof}\, is based on the representation \eqref{eq:1} and is very similar to the proof of Lemma~\ref{lem:1}.

\section{Event-triggering with continuous \nobreak{measurements}}
In Section~\ref{sec:pred} the control signals are sent at $\xi_k=s_k+r_0+\eta_k$, where $r_0+\eta_k$ are sensors-to-controller delays and $s_k$ are measurement sampling instants. In this section we consider the system \eqref{eq:sampSys} without a sensors-to-controller network ($r_0=\eta_k=0$) and with measurements continuously available to the controller. The control law is given by
\begin{equation}\label{cont2}
u(\xi)=Kz(\xi),\quad\xi\ge0,
\end{equation}
where $z(\xi)$ is given by \eqref{pred} with $v(\theta)$ to be defined. To obtain the time instants $\{\xi_k\}$ when a continuously changing control signal $u(\xi)$ is sampled and sent through a controller-to-actuators network, we use a switching approach to event-triggered control \cite{Selivanov2015c}. Namely, we choose $\xi_0=0$,
\begin{equation}\label{evTr}
\xi_{k+1}=\min\{\xi\geq \xi_k+h\,|\,(u(\xi_k)-u(\xi))^T\Omega(u(\xi_k)-u(\xi))\geq\sigma u^T(\xi)\Omega u(\xi)\},
\end{equation}
where a matrix $\Omega\geq0$ and scalars $h>0$, $\sigma\ge0$ are event-triggering parameters. According to \eqref{evTr}, after the controller sends out the control signal $u(\xi_k)$, it waits for at least $h$ seconds. Then it starts to continuously check the event-triggering rule and sends the next control signal when the event-triggering condition is violated. The idea of a switching approach to event-triggered control is to present the closed-loop system as a switching between a system with sampling $h$ and a system with event-triggering mechanism. This allows to ensure large inter-event times and reduce the amount of sent signals \cite{Selivanov2015c}.

Calculating $\dot z$ given by \eqref{pred} in view of \eqref{eq:sampSys} we obtain \eqref{eq:vzSys} (with $r_0=\eta_k=0$). Similar to Section~\ref{sec:evtr} depending on the value of $\mu_M$ one should choose different functions $v(\theta)$.

1. Let $\mu_M=0$. For $v(\theta)$ given in \eqref{vdef}, eq. \eqref{eq:vzSys} takes the form
\begin{equation*}
\dot z(t)=Az(t)+Bu(\xi_k),\quad t\in[\xi_k,\xi_{k+1}).
\end{equation*}
Following \cite{Selivanov2015c} we present the latter system as a switching between two systems:
\begin{equation}\label{eq:7}
\begin{aligned}
\dot z(t)&=Az(t)+BKz(t-\tau_3(t)),&&t\in[\xi_k,\xi_k+h),\\
\dot z(t)&=(A+BK)z(t)+Be_2(t),&&t\in[\xi_k+h,\xi_{k+1}),
\end{aligned}
\end{equation}
where the initial conditions are given by \eqref{inCond}, and
\begin{equation*}
\begin{aligned}
\tau_3(t)&\triangleq t-\xi_k\leq h,&&t\in[\xi_k,\xi_k+h),\\
e_2(t)&\triangleq Kz(\xi_k)-Kz(t),&&t\in[\xi_k+h,\xi_{k+1})
\end{aligned}
\end{equation*}
and \eqref{evTr} implies
\begin{equation*}
0\le\sigma z^T(t)K^T\Omega Kz(t)-e_2^T(t)\Omega e_2(t),\quad t\in[\xi_k+h,\xi_{k+1}).
\end{equation*}

2. Let $\mu_k\not\equiv0$. As it has been explained in Section~\ref{sec:evtr}, in this case it is reasonable not to include the error due to triggering in the definition of $v(\theta)$. Therefore, we take
\begin{equation}\label{v}
v(\xi)\triangleq u(\xi)=Kz(\xi),\quad \xi\ge0.
\end{equation}
Then by calculating $\dot z$ we obtain
\begin{equation}\label{sys1}
\begin{aligned}
\dot z(t)&=Az(t)+BKz(t)-e^{Ar_1}BKz(t-r_1),&&t\in[0,t_0), \\
\dot z(t)&=Az(t)+BKz(t)+e^{Ar_1}BK[z(\xi_k)-z(t-r_1)],&&t\in[t_k,t_{k+1}),\quad k\in\Z_+.
\end{aligned}
\end{equation}
Further analysis of the system \eqref{sys1} is based on a switching approach to event-triggered control~\cite{Selivanov2015c}. Define
\begin{equation*}
t_{-1}^*\triangleq\min\{h,t_0\},\quad t_k^*\triangleq\min\{t_k+h,t_{k+1}\}\text{ for }k\in\Z_+.
\end{equation*}
We have $z(t-r_1-\tau_4(t))=0$ for $t\in[0,t_{-1}^*)$ and $z(\xi_k)=z(t-r_1-\tau_4(t))$ for $t\in[t_k,t_k^*)$, where
\begin{equation*}
\tau_4(t)\triangleq\left\{\begin{aligned}
&\mu_0,&&t\in[0,t_{-1}^*), \\
&t-\xi_k-r_1,&&t\in[t_k,t_k^*).
\end{aligned}\right.
\end{equation*}
Note that $\tau_4(t)\le\tilde\tau\triangleq h+\mu_M$. Further, $Kz(t-r_1-\mu(t))+e_3(t)=0$ for $t\in[t_{-1}^*,t_0)$ and $Kz(\xi_k)=Kz(t-r_1-\mu(t))+e_3(t)$ for $t\in[t_k^*,t_{k+1})$, where
\begin{equation*}
\begin{aligned}
\mu(t)&\triangleq\left\{\begin{aligned}
&\mu_0,&&t\in[t_{-1}^*,t_0),\\
&\mu_k+(t-t_k-h)\frac{\mu_{k+1}-\mu_k}{t_{k+1}-t_k-h},&&t\in[t_k^*,t_{k+1}),\\
\end{aligned}\right.\\
e_3(t)&\triangleq\left\{\begin{aligned}
&0,&&t\in[t_{-1}^*,t_0), \\
&Kz(\xi_k)-Kz(t-r_1-\mu(t)),&&t\in[t_k^*,t_{k+1}).
\end{aligned}\right.
\end{aligned}
\end{equation*}
The function $\mu(t)$ is chosen so that $t-r_1-\mu(t)\in[\xi_k+h,\xi_{k+1})$ for $t\in[t_k^*,t_{k+1})$, therefore, \eqref{evTr} implies
\begin{equation}\label{evtr2}
0\le\sigma z^T(t-r_1-\mu(t))K^T\Omega Kz(t-r_1-\mu(t))-e_3^T(t)\Omega e_3(t)
\end{equation}
for $t\in[t_{k-1}^*,t_k)$ with $k\in\Z_+$.

Finally, the system \eqref{sys1} is presented in the form
\begin{align}\label{sys2}
\dot z(t)&=(A+BK)z(t)+e^{Ar_1}BK[z(t-r_1-\tau_4(t))-z(t-r_1)],\quad t\in[0,t_{-1}^*)\cup[t_k,t_k^*),\\
\label{sys3}\dot z(t)&=(A+BK)z(t)+e^{Ar_1}BK[z(t-r_1-\mu(t))-z(t-r_1)]+e^{Ar_1}Be_3(t),\,t\in[t_{-1}^*,t_0)\cup[t_k^*,t_{k+1})
\end{align}
with \eqref{inCond} and $0\leq\tau_4(t)\leq\tilde\tau=h+\mu_M$, $0\leq\mu(t)\leq\mu_M$.

\begin{lem}\label{lem:2}
	For given $\mu_M\ge0$ and $\alpha>0$ let there exist an $n\times n$ matrix $P>0$, $n\times n$ non-negative matrices $S$, $S_0$, $S_1$, $R_0$, $R_1$, an $m\times m$ matrix $\Omega\ge0$, and $n\times n$ matrices $P_2$, $P_3$, $G_0$, $G_1$ such that
	\begin{equation*}
	\Sigma\leq0,\quad\Xi\le0,\quad\left[\begin{smallmatrix}
	R_0 & G_0 \\
	* & R_0
	\end{smallmatrix}\right]\ge0,\quad
	\left[\begin{smallmatrix}
	R_1 & G_1 \\
	* & R_1
	\end{smallmatrix}\right]\ge0,
	\end{equation*}
	where $\Sigma=\{\Sigma_{ij}\}$ and $\Xi=\{\Xi_{ij}\}$ are the symmetric matrices composed from the matrices
	\begin{align*}
	\Sigma_{11}&=\Xi_{11}=2\alpha P+S+P_2^T(A+BK)+(A+BK)^TP_2,&\Sigma_{56}&=\tilde\rho(R_1-G_1),\\
	\Sigma_{12}&=\Xi_{12}=P-P_2^T+(A+BK)^TP_3,&\Sigma_{66}&=-\tilde\rho(S_1+R_1),\\
	\Sigma_{15}&=\Xi_{14}=-\Sigma_{13}=-\Xi_{13}=P_2^Te^{Ar_1}BK,&\Xi_{17}&=P_2^Te^{Ar_1}B,\\
	\Sigma_{22}&=\Xi_{22}=\mu_M^2R_0+h^2R_1-P_3-P_3^T,&\Xi_{27}&=P_3^Te^{Ar_1}B,\\
	\Sigma_{25}&=\Xi_{24}=-\Sigma_{23}=-\Xi_{23}=P_3^Te^{Ar_1}BK,&\Xi_{34}&=\Xi_{45}=\rho_M(R_0-G_0),\\
	\Sigma_{33}&=\Xi_{33}=e^{-2\alpha r_1}(S_0-S)-\rho_MR_0,&\Xi_{35}&=\rho_MG_0,\\
	\Sigma_{34}&=\rho_MR_0,&\Xi_{44}&=-\Xi_{34}-\Xi_{34}^T+\sigma K^T\Omega K,\\
	\Sigma_{44}&=-\rho_M(R_0+S_0-S_1)-\tilde\rho R_1,&\Xi_{55}&=\rho_M(S_1-S_0-R_0)-\tilde\rho R_1,\\
	\Sigma_{45}&=\tilde\rho(R_1-G_1),&\Xi_{56}&=\tilde\rho R_1,\\
	\Sigma_{46}&=\tilde\rho G_1,&\Xi_{66}&=-\tilde\rho(S_1+R_1),\\	
	\Sigma_{55}&=-\Sigma_{45}-\Sigma_{45}^T,&\Xi_{77}&=-\Omega,\\
	\end{align*}
	$\tilde\rho=e^{-2\alpha(r_1+\tilde\tau)}$, $\rho_M=e^{-2\alpha(r_1+\mu_M)}$, other blocks are zero matrices. Then the system \eqref{inCond}, \eqref{sys2}, \eqref{sys3} with $\xi_k$ given by \eqref{evTr} is exponentially stable with a decay rate~$\alpha$ (i.e. \eqref{zdec} holds).
\end{lem}
{\em Proof} is given in Appendix~\ref{app:lem:2}.

\begin{thm}[Continuous event-triggering]\label{th:2}
	Under the conditions of Lemma~\ref{lem:2} the system \eqref{eq:sampSys}, \eqref{pred}, \eqref{cont2}, \eqref{evTr} with $v(\theta)$ given by \eqref{v} is exponentially stable with a decay rate $\alpha$ (i.e. \eqref{xdec} holds).
\end{thm}
{\em Proof}\, is similar to the proof of Theorem~\ref{th:1}.

\begin{rem}\label{rem:con}
	The control law \eqref{pred}, \eqref{cont2} with $v(\theta)$ given by \eqref{v} requires the knowledge of $z(t)$ for any $t\ge0$. To obtain $z(t)$ during the evolution of the system \eqref{eq:sampSys}, \eqref{pred}, \eqref{evTr}, \eqref{cont2}, \eqref{v} one has to solve the differential equation
	\begin{equation*}
	\begin{aligned}
	\dot z(t)&=(A+BK)z(t)-e^{Ar_1}BKz(t-r_1),&&t\in[0,t_0),\\
	\dot z(t)&=(A+BK)z(t)+e^{Ar_1}BK[z(\xi_k)-z(t-r_1)],&&t\in[t_k,t_{k+1})
	\end{aligned}
	\end{equation*}
	with $z(0)=e^{Ar_1}x(0)$ and $z(t)=0$ for $t<0$.
\end{rem}

\begin{prop}\label{prop:3}
	For $\mu_M=0$ and a given $\alpha>0$, if there exist $n\times n$ matrices $P>0$, $S\ge0$, $R\ge0$, an $m\times m$ matrix $\Omega\ge0$, and $n\times n$ matrices $P_2$, $P_3$, $G$ such that
	\begin{equation*}
	M\leq0,\quad N\leq0,\quad\left[\begin{smallmatrix}
	R & G \\
	* & R
	\end{smallmatrix}\right]\ge0,
	\end{equation*}
	where $M=\{M_{ij}\}$ and $N=\{N_{ij}\}$ are the symmetric matrices composed from the matrices
	\begin{align*}
	M_{11}&=2\alpha P+S-\rho_hR+P_2^TA+A^TP_2,&	M_{33}&=-M_{34}-M_{34}, \\
	N_{11}&=2\alpha P+S-\rho_hR+\sigma K^T\Omega K&M_{44}&=-\rho_h(S+R),\\
	&+P_2^T(A+BK)+(A+BK)^TP_2,&N_{12}&=P-P_2^T+(A+BK)^TP_3,\\
	M_{12}&=P-P_2^T+A^TP_3,&N_{13}&=\rho_hR,\\
	M_{13}&=\rho_h(R-G)+P_2^TBK,&N_{14}&=P_2^TB,\\
	M_{14}&=\rho_hG,&N_{22}&=h^2R-P_3-P_3^T,\\
	M_{22}&=h^2R-P_3-P_3^T,&N_{24}&=P_3^TB,\\
	M_{23}&=P_3^TBK,&N_{33}&=-\rho_h(S+R),\\
	M_{34}&=\rho_h(R-G),&N_{44}&=-\Omega,\\
	\end{align*}
	$\rho_h=e^{-2\alpha h}$, other blocks are zero matrices, then the system \eqref{eq:sampSys}, \eqref{pred}, \eqref{evTr}, \eqref{cont2} with $v(\theta)$ given by \eqref{v} is exponentially stable with a decay rate~$\alpha$.
\end{prop}
{\em Proof}\, is based on the representation \eqref{eq:7} and is very similar to the proof of Lemma~\ref{lem:2}.

\begin{rem}
	MATLAB codes for solving the LMIs of Theorems~\ref{th:1}, \ref{th:2} and Propositions~\ref{prop:1_2}, \ref{prop:3} are available at \url{https://github.com/AntonSelivanov/Aut16a}
\end{rem}

\begin{rem}\label{rem:design}
	Let us set $P_3=\varepsilon_1 P_2$, $\Omega=\varepsilon_2 I_m$ and multiply LMIs of Lemmas~\ref{lem:1},~\ref{lem:2}, Propositions~\ref{prop:1_2}, \ref{prop:3} by $I\otimes P_2^{-1}$ and its transposed from the right and the left, respectively. By denoting $\bar P_2=P_2^{-1}$, $Y=K\bar P_2$ and applying Schur complement to $\sigma Y^T\Omega Y$, we obtain LMIs with tuning parameters $\varepsilon_1$, $\varepsilon_2$ that allow to find controller gain $K=Y\bar P_2^{-1}$. Since requirements $P_3=\varepsilon_1 P_2$, $\Omega=\varepsilon_2 I_m$ may be restrictive, after obtaining $K$ one should use Lemmas~\ref{lem:1},~\ref{lem:2} or Propositions~\ref{prop:1_2}, \ref{prop:3} to obtain larger bound for time-delays and a decay rate. For the details on the LMI-based design see \cite{Fridman2014,Suplin2007}.
\end{rem}

%\begin{rem}
%The system models of this paper allow to consider packet dropouts in communication networks. Let $d^{sc}$ and $d^{ca}$ be the maximum numbers of consecutive packet losses in the sensors-to-controller and controller-to-actuators networks, respectively. Suppose the sensors send out the measurements with the maximum sampling period $h_s$. Denoting by $s_k$ the successfully transmitted measurement we obtain that $h=\max_k(s_{k+1}-s_k)\leq(d^{sc}+1)h_s$.
%
%The lost control signals can be thought of as messages that arrived to the actuator together with later signals and, therefore, were ignored. Namely, if $u(\xi_i)$ and $u(\xi_{i+d})$ were successfully transmitted and $u(\xi_j)$ with $i<j<i+d$ were lost then $\mu_j$ are undefined. We set
%\begin{multline*}
%\hspace{-.3cm}\mu_j\triangleq\xi_{i+d}-\xi_j+\mu_{i+d}=(s_{i+d}+r_0+\eta_{i+d})-(s_j+r_0+\eta_j)+\mu_{i+d}\\
%\leq hd^{ca}+\eta_M+\mu_M\triangleq\tilde{\mu}_M,\quad i<j<i+d
%\end{multline*}
%to obtain
%\begin{equation*}
%t_j=\xi_j+r_0+\mu_j=\xi_{i+d}+r_0+\mu_{i+d}=t_{i+d},\quad i<j<i+d.
%\end{equation*}
%Since $\{t_k\}$ satisfies \eqref{etamu}, the results of this paper can be applied where instead of $\mu_M$ one should use $\tilde{\mu}_M$.
%\end{rem}

\section{Example: inverted pendulum on a cart}
\begin{table*}
	\renewcommand{\arraystretch}{1}
	\centering
	\begin{tabular}{l|c|c|c|c|c|c}
		&\multicolumn{3}{|c|}{$r_0=0.2$, $\eta_M=0.01$} & \multicolumn{3}{|c}{$r_0=\eta_M=0$} \\ \cline{2-7}
		& $\sigma$ & $h$ & SCS & $\sigma$ & $h$ & SCS \\ \hline
		Sampled predictor \eqref{vdef}, \eqref{cont} & $0$ & $0.0369$ & $543$ & $0$ & $0.0646$ & $310$ \\
		Sampled event-triggering \eqref{EvTrRe}, \eqref{EvTrCont} & $0.01$ & $0.0315$ & $116$ & $0.07$ & $0.046$ & $56$ \\
		Continuous predictor \eqref{pred}, \eqref{v} & --- & --- & --- & $0$ & $0.105$ & $191$ \\
		Switching event-triggering \eqref{evTr} & --- & --- & --- & $0.13$ & $0.105$ & $48$
	\end{tabular}
	\caption{Sent control signals (SCS) for different control strategies ($\alpha=0.01$, $r_1=0.2$, $\mu_M=0.01$)}\label{table:sc}
\end{table*}
Following \cite{Wang2009b} we consider an inverted pendulum on a cart controlled through a network described by \eqref{eq:inSys} with
\begin{equation}\label{sysPar}
A=\begin{bmatrix}
0 & 1 & 0 & 0 \\
0 & 0 & -mgM^{-1} & 0 \\
0 & 0 & 0 & 1 \\
0 & 0 & g/l & 0
\end{bmatrix},\quad
B=\begin{bmatrix}
0 \\ M^{-1} \\ 0 \\ -(Ml)^{-1}
\end{bmatrix},
\end{equation}
where $M=10$ kg is the cart mass, $m=1$ kg is the bob mass, $l=3$ m is the arm length and $g=10$ m/s$^2$ is the gravitational acceleration. The state $x=(y,\dot y,\theta,\dot\theta)^T$ is combined of cart's position $y$, cart's velocity $\dot y$, bob's angle $\theta$ and bob's angular velocity $\dot\theta$. For such parameters the open-loop system is unstable and can be stabilized by the control law $u(t)=Kx(t)$ with $K=[2, 12, 378, 210]$. In what follows we compare different control strategies proposed in this paper.

We start by considering a system with both sensors-to-controller and controller-to-actuators networks. The numerical simulations show that the system \eqref{eq:sampSys}, \eqref{sysPar} under the controller $u(t)=Kx(t)$ (without a predictor) is not stable for $r_0=r_1=0.1$, $h=0.0369$, and $\eta_M=\mu_M=0$. The conditions of Corollary~\ref{cor:2} are satisfied for the same $h$ and larger $r_0=r_1=0.2$, $\eta_M=\mu_M=0.01$, whereas the decay rate is $\alpha=0.01$. That is, the predictor-based control admits significantly larger network delays. Furthermore, this implies that within $20$ seconds of simulation $\left\lfloor20/h\right\rfloor+1=543$ signals are sent through each network in the system \eqref{eq:sampSys}, \eqref{sysPar} under the predictor-based controller \eqref{vdef}, \eqref{cont} ($\lfloor\cdot\rfloor$ stands for the integer part). For the system \eqref{eq:startSys}, \eqref{sysPar} under the event-triggered controller \eqref{vdef}, \eqref{cont}, \eqref{EvTrRe}, \eqref{EvTrCont} with $\sigma=0.01$ Theorem~\ref{th:1} gives $h=0.0315$. This bound is smaller than the one given by Corollary~\ref{cor:2}, which means that the event-triggered control requires the measurements $x(s_k)$ to be sent more often but allows to reduce the amount of sent control values $u(\xi_k)$. To obtain the amount of sent signals under the event-triggered control, we perform numerical simulations with $x(0)=[0.98, 0, 0.2, 0]$ and random $\eta_k$, $\mu_k$ satisfying \eqref{etamu}. The results are given in Table~\ref{table:sc}. As one can see event-triggering allows to reduce the workload of the controller-to-actuators network by more than $75\%$. The total amount of signals sent through both sensors-to-controller and controller-to-actuators networks is $543\cdot2=1086$ for the predictor-based controller \eqref{vdef}, \eqref{cont} and $\lfloor20/h\rfloor+1+116=751$ for the event-triggered controller \eqref{vdef}, \eqref{cont}, \eqref{EvTrRe}, \eqref{EvTrCont}.

Now we consider a system with only a controller-to-actuators network ($r_0=\eta_M=0$) and continuous measurements. For this case one can apply sampled predictor-based controller \eqref{vdef}, \eqref{cont} or sampled event-triggered controller \eqref{vdef}, \eqref{cont}, \eqref{EvTrRe}, \eqref{EvTrCont} (with $s_k=\xi_k$). The sampled approach simplifies the calculation of the integral term in \eqref{pred} but does not take advantage of the continuously available measurements. Indeed, as one can see from Table~\ref{table:sc} continuous predictor \eqref{pred}, \eqref{v} without event-triggering ($\sigma=0$ in \eqref{evTr}) reduces the network workload compared to the sampled predictor by almost $40\%$.

To compare the sampled event-triggering mechanism \eqref{vdef}, \eqref{pred}, \eqref{EvTrRe}, \eqref{EvTrCont} and the switching event-triggering mechanism \eqref{pred}, \eqref{evTr}, \eqref{v}, for $\alpha=0.01$ and each value of $\sigma=0.01,0.02,\ldots,1$ we apply Theorems~\ref{th:1} and \ref{th:2} to find the maximum allowable $h$. Then we perform numerical simulations for each pair of $(\sigma,h)$ with $\mu_k$ subject to \eqref{etamu} ($r_1=0.2$, $\mu_M=0.01$) and choose the pair $(\sigma,h)$ that leads to the smallest amount of sent control signals. In Table~\ref{table:sc} one can see that both event-triggering mechanisms significantly reduce the amount of sent control signals. The switching event-triggering reduces the network workload by almost $15\%$ compared to the sampled event-triggering.

\section{Conclusions}
We considered predictor-based control of NCSs with uncertain network delays. For the event-triggered control we showed that one should use different predictor models depending on the value of the controller-to-actuators delay uncertainty. To take advantage of the continuously available measurements in the systems with only a controller-to-actuators network, we considered a continuous-time predictor with a switching event-triggering mechanism. For the proposed control strategies we obtained LMI-based stability conditions that guaranty the desired exponential decay rate of convergence and allow to find appropriate controller gains. An example of inverted pendulum on a cart demonstrates that event-triggering mechanism allows to reduce the network workload and in those cases where the continuous-time predictor can be applied it has some advantages over the sampled one.

\bibliographystyle{IEEEtran}
\bibliography{library}
\appendix
\section{Proof of Lemma~\ref{lem:1}}\label{app:lem:1}
For $t\geq\tau_M$ we have
\begin{equation}\label{eq1}
\begin{aligned}
\dot V_P+2\alpha V_P&=2z^T(t)P\dot z(t)+2\alpha z^T(t)Pz(t),\\
\dot V_{S_0}+2\alpha V_{S_0}&=z^T(t)S_0z(t)-e^{-2\alpha\bar\tau}z^T(t-\bar\tau)S_0z(t-\bar\tau), \\
\dot V_S+2\alpha V_S&=e^{-2\alpha\bar\tau}z^T(t-\bar\tau)Sz(t-\bar\tau)-e^{-2\alpha(r_0+r_1)}z^T(t-r_0-r_1)Sz(t-r_0-r_1), \\
\dot V_{S_1}+2\alpha V_{S_1}&=e^{-2\alpha(r_0+r_1)}z^T(t-r_0-r_1)S_1z(t-r_0-r_1)-e^{-2\alpha \tau_M}z^T(t-\tau_M)S_1z(t-\tau_M).
\end{aligned}
\end{equation}
Using Jensen's inequality \cite{Gu2003}, Park's theorem \cite{Park2011} and taking into account that $\tau_1(t)\le\tau_2(t)$ \cite{Liu2012a} we obtain
\begin{multline}\label{eq2}
\dot V_{R_0}+2\alpha V_{R_0}=\bar\tau^2\dot z^T(t)R_0\dot z(t)-\bar\tau\int_{t-\bar\tau}^te^{2\alpha(s-t)}\dot z^T(s)R_0\dot z(s)\,ds\\
\le\bar\tau^2\dot z^T(t)R_0\dot z(t)-e^{-2\alpha\bar\tau}\left[\begin{smallmatrix}z(t)-z(t-\tau(t))\\
z(t-\tau(t))-z(t-\bar\tau)
\end{smallmatrix}\right]^T
\left[\begin{smallmatrix}R_0 & G_0 \\ G_0^T & R_0
\end{smallmatrix}\right]
\left[\begin{smallmatrix}z(t)-z(t-\tau(t))\\
z(t-\tau(t))-z(t-\bar\tau)
\end{smallmatrix}\right],
\end{multline}
\begin{multline}\label{eq3}
\dot V_{R_1}+2\alpha V_{R_1}=(\tau_M-r_0-r_1)^2\dot z^T(t)R_1\dot z(t)-(\tau_M-r_0-r_1)\int_{t-\tau_M}^{t-r_0-r_1}e^{2\alpha(s-t)}\dot z^T(s)R_1\dot z(s)\,ds\\
\le
(\tau_M-r_0-r_1)^2\dot z^T(t)R_1\dot z(t)-e^{-2\alpha\tau_M}\left[\begin{smallmatrix}
z(t-r_0-r_1)-z(t-\tau_1(t))\\
z(t-\tau_1(t))-z(t-\tau_2(t))\\
z(t-\tau_2(t))-z(t-\tau_M)
\end{smallmatrix}\right]^T
\left[\begin{smallmatrix}
R_1 & G_1 & G_2 \\
* & R_1 & G_3 \\
* & * & R_1
\end{smallmatrix}\right]
\left[\begin{smallmatrix}
z(t-r_0-r_1)-z(t-\tau_1(t))\\
z(t-\tau_1(t))-z(t-\tau_2(t))\\
z(t-\tau_2(t))-z(t-\tau_M)
\end{smallmatrix}\right].
\end{multline}
We use the following descriptor representation of \eqref{eq:zSys}
\begin{multline}\label{eq4}
0=2[z^T(t)P_2^T+\dot z^T(t)P_3^T]\bigl[-\dot z(t)+Az+BKz(t-\tau(t))\\
+e^{A(r_0+r_1)}B\bigl(e_1(t)+Kz(t-\tau_2(t))-Kz(t-\tau_1(t))\bigr)\bigr].
\end{multline}
By summing up \eqref{EvTrCompensation}, \eqref{eq1}--\eqref{eq4} we obtain
\begin{equation*}
\dot V+2\alpha V\leq\varphi^T\Phi\varphi\leq0,
\end{equation*}
where $\varphi=\operatorname{col}\{z(t), \dot z(t), z(t-\tau(t)), z(t-\bar\tau), z(t-r_0-r_1), z(t-\tau_1(t)), z(t-\tau_2(t)), z(t-\tau_M), e_1(t)\}$. This implies $\dot V\leq-2\alpha V$ and, therefore,
\begin{equation}\label{eq:8}
V(t)\leq e^{-2\alpha(t-\tau_M)}V(\tau_M), \quad t\ge\tau_M.
\end{equation}
Define $z_t=z(t+\theta)$, $\theta\in[-\tau_M,0]$ and $\|z_t\|_{PC}=\max_{\theta\in[-\tau_M,0]}|z(t+\theta)|$.
For $t\ge0$ function $\|z_t\|_{PC}$ is continuous in $t$ and \eqref{eq:zSys}, \eqref{EvTrCompensation} imply $|\dot z(t)|\le m\|z_t\|_{PC}$ for some $m>0$. Therefore,
\begin{equation*}
\begin{array}{c}
\|z_t\|_{PC}\leq|z(0)|+\int_0^tm\|z_s\|_{PC}\,ds,\quad t\ge0.
\end{array}
\end{equation*}
By the Gronwall-Bellman Lemma this implies
\begin{equation}\label{eq:10}
\|z_t\|_{PC}\leq|z(0)|e^{mt},\quad t\ge0.
\end{equation}
Since $|\dot z(t)|\le m\|z_t\|_{PC}$, there exists $c_1$ such that $V(\tau_M)\leq c_1\|z_{\tau_M}\|_{PC}^2\leq c_1|z(0)|^2e^{2m\tau_M}$. Since $|z(t)|^2\lambda_{\min}(P)\leq V(t)$, \eqref{eq:8} and \eqref{eq:10} imply \eqref{zdec} for some $M>0$.

\section{Proof of Theorem~\ref{th:1}}\label{app:th:1}
From \eqref{vdef}, \eqref{pred}, \eqref{eq:9} we have
\begin{equation*}
\begin{array}{l}
x(t)=e^{-A(r_0+r_1)}z(t)-\int_{t-r_1}^{t+r_0}e^{A(t-r_1-\theta)}BKz(\theta-r_0-\tau(\theta-r_0))\,d\theta,\quad t\ge0,
\end{array}
\end{equation*}
where $z$ satisfies \eqref{inCond}, \eqref{eq:zSys}. By Lemma~\ref{lem:1} \eqref{zdec} holds, thus
\begin{equation*}
\begin{array}{c}
|x(t)|\leq Ce^{-\alpha t}|z(0)|=Ce^{-\alpha t}\left\|e^{A(r_0+r_1)}\right\||x(0)|.
\end{array}
\end{equation*}

\section{Proof of Lemma~\ref{lem:2}}\label{app:lem:2}
For $t\ge r_1+\tilde\tau$ ($\tilde\tau=h+\mu_M$) consider the functional
\begin{equation*}
V=V_P+V_S+V_{S_0}+V_{R_0}+V_{S_1}+V_{R_1},
\end{equation*}
%where
\begin{equation*}
\begin{array}{l}
V_P=z^T(t)Pz(t),\\
V_S=\int_{t-r_1}^te^{2\alpha(s-t)}z^T(s)Sz(s)\,ds,\\
V_{S_0}=\int_{t-r_1-\mu_M}^{t-r_1}e^{2\alpha(s-t)}z^T(s)S_0z(s)\,ds,\\
V_{R_0}=\mu_M\int_{-r_1-\mu_M}^{-r_1}\int_{t+\theta}^te^{2\alpha(s-t)}\dot z^T(s)R_0\dot z(s)\,ds\,d\theta,\\
V_{S_1}=\int_{t-r_1-\tilde\tau}^{t-r_1-\mu_M}e^{2\alpha(s-t)}z^T(s)S_1z(s)\,ds,\\
V_{R_1}=h\int_{-r_1-\tilde\tau}^{-r_1-\mu_M}\int_{t+\theta}^te^{2\alpha(s-t)}\dot z^T(s)R_1\dot z(s)\,ds\,d\theta.
\end{array}
\end{equation*}
We have
\begin{equation}\label{eq:2}
\begin{array}{l}
\dot V_P+2\alpha V_P=2z^T(t)P\dot z(t)+2\alpha z^T(t)Pz(t),\\
\dot V_S+2\alpha V_S=z^T(t)Sz(t)-e^{-2\alpha r_1}z^T(t-r_1)Sz(t-r_1), \\
\dot V_{S_0}+2\alpha V_{S_0}=e^{-2\alpha r_1}z^T(t-r_1)S_0z(t-r_1)-e^{-2\alpha(r_1+\mu_M)}z^T(t-r_1-\mu_M)S_0z(t-r_1-\mu_M), \\
\dot V_{S_1}+2\alpha V_{S_1}=e^{-2\alpha(r_1+\mu_M)}z^T(t-r_1-\mu_M)S_1\times\\
\hspace{2.8cm}z(t-r_1-\mu_M)-e^{-2\alpha(r_1+\tilde\tau)}z^T(t-r_1-\tilde\tau)S_1z(t-r_1-\tilde\tau), \\
\dot V_{R_0}+2\alpha V_{R_0}=\mu_M^2\dot z^T(t)R_0\dot z(t)-\mu_M\int_{t-r_1-\mu_M}^{t-r_1}e^{2\alpha(s-t)}\dot z^T(s)R_0\dot z(s)\,ds, \\
\dot V_{R_1}+2\alpha V_{R_1}=h^2\dot z^T(t)R_1\dot z(t)-h\int_{t-r_1-\tilde\tau}^{t-r_1-\mu_M}e^{2\alpha(s-t)}\dot z^T(s)R_1\dot z(s)\,ds.
\end{array}
\end{equation}
I. For $t\in[t_k^*,t_{k+1})$ we have
\begin{equation}\label{eq6}
\begin{aligned}
0=2[z^T(t)P_2^T+\dot z^T(t)P_3^T][-\dot z(t)+(A+BK)z(t)
+e^{Ar_1}B(Kz(t-r_1-\mu(t))-Kz(t-r_1)+e_3(t))].
\end{aligned}
\end{equation}
To compensate the term $z(t-r_1-\mu(t))$ using Jensen's inequality and Park's theorem we derive
\begin{multline}\label{eq:5}
-\mu_M\int_{t-r_1-\mu_M}^{t-r_1}e^{2\alpha(s-t)}\dot z^T(s)R_0\dot z(s)\,ds\le e^{-2\alpha(r_1+\mu_M)}\times\\
\left[\begin{smallmatrix}
z(t-r_1)-z(t-r_1-\mu(t)) \\
z(t-r_1-\mu(t))-z(t-r_1-\mu_M)
\end{smallmatrix}\right]^T\left[\begin{smallmatrix}
R_0 & G_0 \\ G_0^T & R_0
\end{smallmatrix}\right]\left[\begin{smallmatrix}
z(t-r_1)-z(t-r_1-\mu(t)) \\
z(t-r_1-\mu(t))-z(t-r_1-\mu_M)
\end{smallmatrix}\right],
\end{multline}
\begin{multline}\label{eq:6}
-\mu_M\int_{t-r_1-\mu_M}^{t-r_1}e^{2\alpha(s-t)}\dot z^T(s)R_0\dot z(s)\,ds \\ \leq e^{-2\alpha(r_1+\tilde\tau)}\left[z(t-r_1-\mu_M)-z(t-r_1-\tilde\tau)\right]^TR_1[z(t-r_1-\mu_M)-z(t-r_1-\tilde\tau)].
\end{multline}
By summing up \eqref{evtr2}, \eqref{eq:2}, \eqref{eq6} in view of \eqref{eq:5}, \eqref{eq:6} we obtain $\dot V+\alpha V\leq\xi^T\Xi\xi\le0$,
where $\xi=\operatorname{col}\{z(t), \dot z(t), z(t-r_1), z(t-r_1-\mu(t)), z(t-r_1-\mu_M), z(t-r_1-\tilde\tau), e_3(t)\}$.

For $t\in[t_k,t_k^*)$ the system \eqref{sys2} with $\tau_4(t)\in[0,\mu_M)$ is described by \eqref{sys3} with $e_3(t)=0$ satisfying \eqref{evtr2}.

II. For $t\in[t_k,t_k^*)$, $\tau_4(t)\in[\mu_M,\mu_M+h)$ we have
\begin{multline}\label{eq5}
0=2\left[z^T(t)P_2^T+\dot z^T(t)P_3^T\right]\left[-\dot z(t)+(A+BK)z(t)\right.\\
\left.+e^{Ar_1}BKz(t-r_1-\tau_4(t))-e^{Ar_1}BKz(t-r_1)\right].
\end{multline}
To compensate the term $z(t-r_1-\tau_4(t))$ using Jensen's inequality and Park's theorem we derive
\begin{multline}\label{eq:3}
-\mu_M\int_{t-r_1-\mu_M}^{t-r_1}e^{2\alpha(s-t)}\dot z^T(s)R_0\dot z(s)\,ds \\ \leq-e^{-2\alpha(r_1+\mu_M)}\left[z(t-r_1)-z(t-r_1-\mu_M)\right]^TR_0[z(t-r_1)-z(t-r_1-\mu_M)],
\end{multline}
\begin{multline}\label{eq:4}
-h\int_{t-r_1-\tilde\tau}^{t-r_1-\mu_M}e^{2\alpha(s-t)}\dot z^T(s)R_1\dot z(s)\,ds\le-e^{-2\alpha(r_1+\tilde\tau)}\times\\
\left[\begin{smallmatrix}
z(t-r_1-\mu_M)-z(t-r_1-\tau_4(t)) \\
z(t-r_1-\tau_4(t))-z(t-r_1-\tilde\tau)
\end{smallmatrix}\right]^T\left[\begin{smallmatrix}
R_1 & G_1 \\ G_1^T & R_1
\end{smallmatrix}\right]\left[\begin{smallmatrix}
z(t-r_1-\mu_M)-z(t-r_1-\tau_4(t)) \\
z(t-r_1-\tau_4(t))-z(t-r_1-\tilde\tau)
\end{smallmatrix}\right].
\end{multline}
By summing up \eqref{eq:2}, \eqref{eq5} in view of \eqref{eq:3}, \eqref{eq:4} we obtain $\dot V+2\alpha V\leq\eta^T\Sigma\eta\le0$, where $\eta=\operatorname{col}\{z(t), \dot z(t), z(t-r_1), z(t-r_1-\mu_M), z(t-r_1-\tau_4(t)), z(t-r_1-\tilde\tau)\}$.

Therefore, we obtain $\dot V\leq-2\alpha V$ for $t\ge r_1+\tilde\tau$. The end of the proof is similar to that of Lemma~\ref{lem:1}.
\end{document}